\documentclass[12pt,a4paper]{article}% For Computer Modern Fonts, or,
\usepackage{latexsym,amssymb}
\usepackage{graphics,times,mathptm}
\usepackage{epsfig}
\def\beq{\begin{equation}}
\def\eeq{\end{equation}}
\def\bea{\begin{eqnarray}}
\def\eea{\end{eqnarray}}
\def\beas{\begin{eqnarray*}}
\def\eeas{\end{eqnarray*}}

\newcommand{\nms}{\normalsize}
\newenvironment{keywords}{ \noindent {\small\bf Key Words}:}{ }

\begin{document}
%\begin{article}

\begin{center}

 {TECHNICAL NOTE}

\vspace{15mm}

 { Efficient Partition of $N$-Dimensional Intervals in the
Framework of One-Point-Based Algorithms}\footnotemark

\vspace{10mm}

\footnotetext{\nms This research was  supported by  the Italian
Fund of Fundamental Research (grants FIRB RBNE01WBBB and FIRB
RBAU01JYPN) and by the Russian Fund of Fundamental Research (grant
RFBR 01--01--00587).\vspace{5mm} }

%% Author name

     Ya.D. Sergeyev\footnotemark
\footnotetext{\nms {\it Distinguished Professor}, Dipartimento di
Elettronica, Informatica e Sistemistica, Universit\`a della
Calabria, Via P.~Bucci, Cubo~41-C, 87030 Rende, Cosenza, Italy and
        {\it part-time   Professor}, University of Nizhni Novgorod,
        Nizhni Novgorod, Russia.}

        Communicated by G. Di Pillo

\end{center}

\newpage

\begin{abstract}
\nms In this paper,  the problem of the minimal description  of
the structure of a vector function  $f(x)$  over an
$N$-dimensional interval is   studied. Methods adaptively
subdividing the original interval in smaller subintervals and
evaluating $f(x)$ at only one point within each subinterval are
considered. Two partition strategies traditionally used for
solving this problem are analyzed. A new partition strategy based
on an efficient technique developed for diagonal algorithms is
proposed and studied.
 \end{abstract}

\vspace{2cm}


\newpage

\section{Introduction and Analysis of Traditional Partition\\ Strategies}
The problem of the minimal description of   the behavior of a
vector function
 \beq
  f(x) = (f_1(x), f_2 (x), ..., f_r (x)),   \hspace{5mm}  r \ge 1,  \label{1}
  \eeq
over a  hyperinterval $D \subset \mathbb{R}^N$ can  be  stated in
 various ways  under  different       assumptions regarding
  $f(x)$. The term `minimal description' means that
we want to obtain a knowledge about $f(x)$ by evaluating it in a
minimal  number of trial  points $x \in D$. The problem has a
number of important applications in numerous fields of mathematics
such as optimization,  number theory, numerical integration,
geometric partitioning, and structural description (see Ref.~1 for
discussion and references). Usually, in real-life applications,
the operation of evaluating $f(x)$ requires much time and to
obtain an acceptable solution of the minimal description problem
it is necessary to execute a high number of such evaluations.

Numerous   iterative    processes     proposed in literature (see
Refs.~1--19, etc.) for solving this problem can be distinguished
in dependence on the way they combine the following four features:
(i)~the strategy used for partition of the region $D$; (ii)~the
way to choose an element (or elements) for the next partition;
(iii)~the number of points   $f(x)$ should be evaluated at the new
subregions obtained after partition; (iv)~the location of these
points within each of the new subregions.

Let us determine the place of our study with respect to the
features (i)--(iv). First,    one-point-based, diagonal,
simplicial, etc. algorithms (see, for example, Refs.~1--19) can be
distinguished  relatively to the feature (i) and (iii).
one-point-based algorithms subsequently subdivide the region $D$
in smaller hyperintervals and evaluate $f(x)$ at one point within
each sub-interval (the terms `cell' and `box' will be also used).
Diagonal methods do the same but evaluate $f(x)$ at two vertices
of each box. Simplicial algorithms partition the region in
simplexes and evaluate the objective vector function  at all their
vertices. The current state-of-art in the field (see Refs.~2, 7,
11, 15) does not allow us to say which of the approaches is worse
or better for a given class of functions.

This paper deals with popular one-point-based algorithms that have
been extensively studied theoretically (mainly from the feature
(ii) point of view (see  Refs.~6, 8, 11, 12, 14, 15, etc.)). They
have also been successfully  applied to solve numerous real-life
problems. For example, interval analysis methods use mainly this
strategy in their work (see, for instance, Refs.~7, 10, 11).
Another important example of their usage comes from the DIRECT
optimization method introduced in Ref.~12 that also has been
employed in a number of industrial applications (see Refs.~5, 17,
18 ).

Peculiarity of this paper consists of the following: it does not
discusses the feature (ii). In contrast, its goal is to show (as
it has been already done for the diagonal methods in Ref.~1) that
partition strategies themselves, independently of the feature
(ii), can influence significantly the number of function
evaluations made by an algorithm. Thus, we concentrate our
attention on the features (i) and (iv) in the framework of
one-point-based algorithms.

Let us analyze two partition strategies traditionally  used in the
one-point-based algorithms. In the first of them, the region $D$
is subdivided in a number of sub-intervals and $f(x)$ is evaluated
at an internal (very often central) point of each of the new
sub-boxes. Then a new sub-cell of $D$ is chosen for partitioning
and the process is repeated. This simple and widely used strategy
has the following drawback. If $f(x)$ has been evaluated at a
point $x^i$ within an interval $D(i)$ only the interval $D(i)$
uses the information obtained from evaluation of $f(x^i)$.

The second traditional partition strategy overcomes this
difficulty. An internal point (indicated in Fig.~1a by the number
$1$) is chosen within the region $D$ that is subdivided in $2^N$
sub-boxes by the hyperplanes orthogonal to coordinate axes and
passing       through this point. In this strategy, the
information obtained at every point is used by $2^N$ intervals.
Unfortunately, such a huge number of sub-boxes creates problems
during managing the description information when the dimension of
the problem, $N$, increases. This strategy has also the second
problem illustrated in Fig.~1b.

Suppose that the interval shown by grey color in Fig.~1a has been
chosen for the next subdivision. It can be seen from Fig.~1b that
the partition executed within this interval at the point $2$
creates  redundancy in the following sense. The one-point-based
algorithms use only {\it one} point for description of $f(x)$ over
each interval $D(i)$. In spite of this, the interval shown by grey
color has {\it two} points where $f(x)$ has been evaluated. One of
them is redundant. In general, usage of this partition strategy
leads to the following result describing the level of the obtained
redundancy.

\noindent {\bf Proposition 1.1} Every partition made by using the
second strategy leads to creation of one or two sub-cells having
$f(x)$ evaluated at two vertices.

\noindent {\bf Proof:} Let us consider the situation where one
interval having $f(x)$  evaluated at two vertices is generated.
This happens when an interval having $f(x)$  evaluated at one of
its vertices (see Fig.~1a) is subdivided. Two intervals having
redundant points are generated when an interval with $f(x)$
evaluated at two of its vertices (see the interval shown in grey
in Fig.~1b) is subdivided.  \rule{5pt}{5pt}

 The analysis given
above shows that a desirable partition strategy should not
generate too many sub-cells during every partition and should be
able to avoid redundant evaluations of the vector function $f(x)$
giving in the same moment to several intervals a possibility to
use the information obtained from every single evaluation of
$f(x)$.

\section{Strategy}
In this paper, a new partition strategy that can be used by
one-point-based algorithms is proposed. It is based on an
efficient partition strategy introduced recently in Ref.~1 for
solving the minimal description problem in the framework of
diagonal algorithms (see Ref.~15). These methods evaluate the
vector function  $f(x)$ at two vertices $a(i)$ and $b(i)$ of each
sub-box $D(i)$ where
\[
D(i) = [a(i),b(i)] = \{ x: a(i) \le x \le b(i) \}.
\]
A high practical efficiency of the new strategy applied for
solving global optimization problems has been shown in Ref.~13.

In order to proceed, let us describe this partition strategy
developed for diagonal algorithms. A cell $D(t)= [a(t),b(t)]$
chosen for subdivision among $L(k)$ cells existing during an
iteration $k$ is split into three equal sub-intervals by two
hyperplanes orthogonal to the longest edge parallel to the $i$th
coordinate axis and passing through the points
 \beq
   u = ( a_1(t),a_2(t),\ldots,a_{i-1}(t), a_i(t)+2(b_i(t)-a_i(t))/3, a_{i+1}(t), \ldots,a_N(t)
   ), \label{2}
 \eeq
 \beq
   v = (b_1(t),b_2(t),\ldots,b_{i-1}(t), b_i(t)+2(a_i(t)-b_i(t))/3,
b_{i+1}(t),\ldots,b_N(t)).   \label{3}
 \eeq
 The cell $D(t)$ is substituted by the new cells $D(t(k)),  D(L(k)+1),$
and      $ D(L(k)+2)$   determined by their vertices
 \[
\hspace*{1cm}  \begin{array}{llr}
      a(t(k)) = a(L(k)+2) = u,    &  b(t(k)) = b(L(k)+1) = v,& \hspace*{2cm}(4a) \\
      a(L(k)+1) = a(t(k)),        &  b(L(k)+1) = v,      & \hspace*{2cm}(4b)     \\
      a(L(k)+2) = u,              &  b(L(k)+2) = b(t(k)). & \hspace*{2cm}(4c)
     \end{array}
 \]
The function  $f(x)$ is (eventually) evaluated at the points $u$
and $v$.

It has been shown in  Ref.~1 that this strategy generates regular
trial meshes in such a way that every cell has exactly two
vertices where the function  $f(x)$ is evaluated. The introduced
regularization allows  to establish links between sub-cells
generated   during different iterations
 eliminating  in this way     possibility of the
redundant trials generation and storage of the related
information.

While using this strategy it becomes possible (see Ref.~1) to
reestablish  information about vicinity of the cells generated
during different iterations and, as the result, to eliminate
redundant storage of the points $a(j), b(j)$ and  results of
evaluations  of the function  $f(x)$ at these       points.
Particularly, it is shown that every vertex where $f(x)$ is
evaluated can belong to different (up to $2^N$) cells. When we
split an interval $D(t)$, we calculate the coordinates of the
vertices corresponding to the three new sub-cells. Particularly,
we are interested   in the vertices $u$ and $v$, the only two
vertices for the three cells from (4) where the function  $f(x)$
should be evaluated.

Instead of an immediate evaluation of the values $f(u)$ and
$f(v)$, we verify the existence of these in the data base because
$f(u)$ and/or $f(v)$ could have been already evaluated during
previous iterations. It is important to mention that two
numerations  (one for the boxes and another for the vertices where
the function  $f(x)$ has been evaluated) have been developed
theoretically in Ref.~1 and successfully applied practically in
Ref.~13. These numerations allow us to calculate the addresses of
$u$ and $v$ in the data base of the vertices from the number of
the box $D(t)$ providing so a direct fast access to the values
$f(u)$ and $f(v)$.

If  both of them have been already evaluated, we simply read
these values from the data base. If only one of them has been
evaluated, we read this value and create a new element in the
data    base for the absent (say $u$) point, evaluate $f(u)$ and
record it       in the element created. In the last case -- both
values are absent -- these operations are   executed for both
points.

Thus, the description information is evaluated at every  vertex
only  once    and then we simply  read  it up to $2^N$ times
instead of  evaluating $f(x)$ and saving the result of this
evaluation and coordinates of the $2^N$ trial points $2^N$ times.

Surprisingly, this strategy developed for the diagonal methods can
be successfully applied for the one-point-based algorithms too.
Instead of evaluating $f(x)$ at two vertices, $u$ and $v$, it is
proposed to do this initially for the vertex $a$ of the region $D$
and then at the vertex $u$ during every splitting (the point $v$
is used just for partitioning goals). Then, the operation of the
verification whether the function  $f(x)$ has been already
evaluated at this point is made by using the fast procedure
developed in Ref.~1.

In order to illustrate the new strategy, let us consider an
example presented in Fig.~2. The first evaluation of the function
$f(x)$ is executed at the vertex $a$ (see Fig.~2a). The second
evaluation is made at the point $2$ being the point $u$ from
(\ref{2}). It can be seen from this figure that we have three
sub-cells having exactly one vertex where $f(x)$ has been
evaluated. Suppose that the interval shown in grey in Fig.~2a has
been chosen for the next splitting (see Fig.~2b). The only
evaluation made at the point $3$ gives us three new intervals. It
seems that the interval chosen for the next partitioning (and
shown again in grey) has a redundant point, i.e., the point $3$.
In reality, this point is not a redundant one but is kept in the
data base and can be used in the future. Fig.~2c shows that the
fourth point coincides with the point $3$, thus, no evaluations
are made, the description information is read from the data base
and we obtain three new smaller intervals gratis. Finally, Fig.~2d
illustrates situation after eleven iterations. It can be seen from
this figure that $21$ intervals have been generated and the
function  has been evaluated at seven points only.

\section{Conclusions}
In this paper, the problem of the minimal description of a
function  $f(x)$ over a  hyperinterval      $D$  has been
considered in the framework of the one-point-based algorithms. An
analysis of the traditional partitioning strategies used by these
methods has been made. It has been shown that the meshes of the
trial points generated by these strategies can generate redundant
points and intervals. Such a redundancy   may lead to a
significant increase of information to be stored in the computer
memory and to the slowing down the description procedure.

A new partition strategy  has been introduced to overcome these
difficulties. It generates regular trial       meshes successfully
responding to the  requirements  of the   minimal description.
During every partitioning it generates three intervals
independently on the problem dimension. Since every point where
the function  has been evaluated can belong up to $2^N$ intervals
and it is possible to establish links between sub-cells generated
during       different iterations, the function  is evaluated at
every trial point only once and then the result of this evaluation
is simply read from the data base many times. As a rule, this fact
very often allows to obtain new partitions without any evaluation
of the function.

\vspace{1cm}

\noindent {\bf \Large References}

\vspace{1cm}

\noindent 1.  \textsc{Sergeyev, Ya. D.},  {\it An Efficient
Strategy for Adaptive Partition of N-Dimensional Intervals in the
Framework of Diagonal Algorithms}, Journal of Optimization Theory
and Applications, Vol. 107, pp. 145--168, 2000.

\noindent 2. \textsc{Agarwal,   P. K.}, {\it Geometric
Partitioning and Its
 Applications},  DIMACS Series in Discrete Mathematics and Theoretical Computer Science, Edited by J. E. Goodman, R. Pollack, and J.
O'Rourke, Vol. 6, pp. 1--37, 1991.

\noindent 3. \textsc{Baritompa,  W. P.,   and  Viitanen, S.}, {\it
 PMB-Parallel         Multidimensional Bisection},
Report 101, Department of Mathematics and Statistics, University
of Canterbury, Christchurch, New Zealand,  1993.

\noindent 4. \textsc{Breiman,  L., and  Cutler, A.}, {\it  A
Deterministic Algorithm for Global Optimization}, Mathematical
Programming, Vol. 58, pp. 179--199, 1993.

\noindent 5. \textsc{Cox, S. E., Haftka, R. T., Baker, C.,
Grossman, B., Mason, W. H., and  Watson, L. T.}, {\it A Comparison
of Global Optimization Methods for the Design of a High-Speed
Civil Transport},  Journal of Global Optimization, Vol. 21, pp.
415--433, 2001.

\noindent 6. \textsc{Csendes,  T. and  Ratz, D.}, {\it Subdivision
Direction Selection in Interval Methods for Global Optimization},
SIAM Journal of Numerical Analysis, Vol. 34, pp. 922--938, 1997.

\noindent 7.   \textsc{Floudas, C., and  Pardalos, P. M.}, {\it
State of the Art in Global Optimization: Computational Methods and
Applications}, Kluwer Academic Publishers, Dordrecht, the
Netherlands, 1995.

\noindent 8. \textsc{Galperin, E. A.}, {\it The Cubic Algorithm,}
 Journal of Mathematical Analysis and Applications, Vol. 112, pp.   635--640, 1985.

\noindent 9. \textsc{Gergel,  V.P.}, {\it A Global Optimization
Algorithm for Multivariate Functions with Lipschitzian First
Derivatives}, Journal of Global Optimization, Vol. 10, pp.
257-281,  1997.

\noindent 10. \textsc{Hansen, E. R.}, {\it Global Optimization
Using Interval Analysis}, Marcel Dekker,   New York, NY, 1992.

\noindent 11. \textsc{Horst, R., and Pardalos, P. M.},
 {\it Handbook of Global Optimization}, Kluwer Academic Publishers, Dordrecht, the
Netherlands, 1995.

\noindent 12. \textsc{Jones, D. R., Perttunen, C. D., and
Stuckman, B. E.}, {\it
 Lipschitzian Optimization without the Lipschitz Constant}, Journal of Optimization
 Theory and Applications, Vol. 79, pp. 157--181, 1993.

\noindent 13. \textsc{Kvasov, D. E., and  Sergeyev, Ya. D.}, {\it
Multidimensional Global Optimization Algorithm Based on Adaptive
Diagonal Curves}, Computational Mathematics and Mathematical
Physics, Vol. 43, pp. 40--56, 2003.

\noindent 14. \textsc{Meewella,  C. C.,  and  Mayne, D. Q.}, {\it
Efficient Domain Partitioning Algorithms for Global Optimization
of Rational and Lipschitz Continuous Functions},  Journal of
  Optimization Theory and Applications, Vol. 61, pp. 247--270, 1989.

\noindent 15. \textsc{Pint\'{e}r, J. D.}, {\it Global Optimization
in Action}, Kluwer Academic Publishers, Dordrecht, the
Netherlands, 1996.

\noindent 16. \textsc{Strongin,  R. G.,  and Sergeyev, Ya. D.},
{\it Global Optimization with Non-Convex Constraints: Sequential
and Parallel Algorithms}, Kluwer Academic Publishers, Dordrecht,
the Netherlands,
 2000.

\noindent 17. \textsc{Verstak, A.,  He, J., Watson, L. T,
Ramakrishnan, N., Shaffer, C. A., Rappaport, T. S., Anderson, C.
R., Bae, K., Jiang, J., and  Tranter, W. H.}, \textit{S4W:
Globally Optimized Design of Wireless Communication Systems},
Proceedings of the Next Generation Software Workshop: 16th
International Parallel and Distributed Processing Symposium, IEEE
Computer Society Press, Fort Lauderdale, Florida, pp. 173--180,
2002.

\noindent 18. \textsc{Watson, L. T., and  Baker, C. A.}, \textit{A
Fully-Distributed Parallel Global Search Algorithm}, Engineering
Computations, Vol. 18, pp. 155--169, 2001.

\noindent 19. \textsc{Wood, G. R.},  {\it The Bisection Method in
Higher Dimensions},  Mathematical Programming, Vol. 55, pp.
319--337, 1992.

 \newpage

\vspace*{5cm}
\begin{figure}[t]
\centerline{\psfig{figure=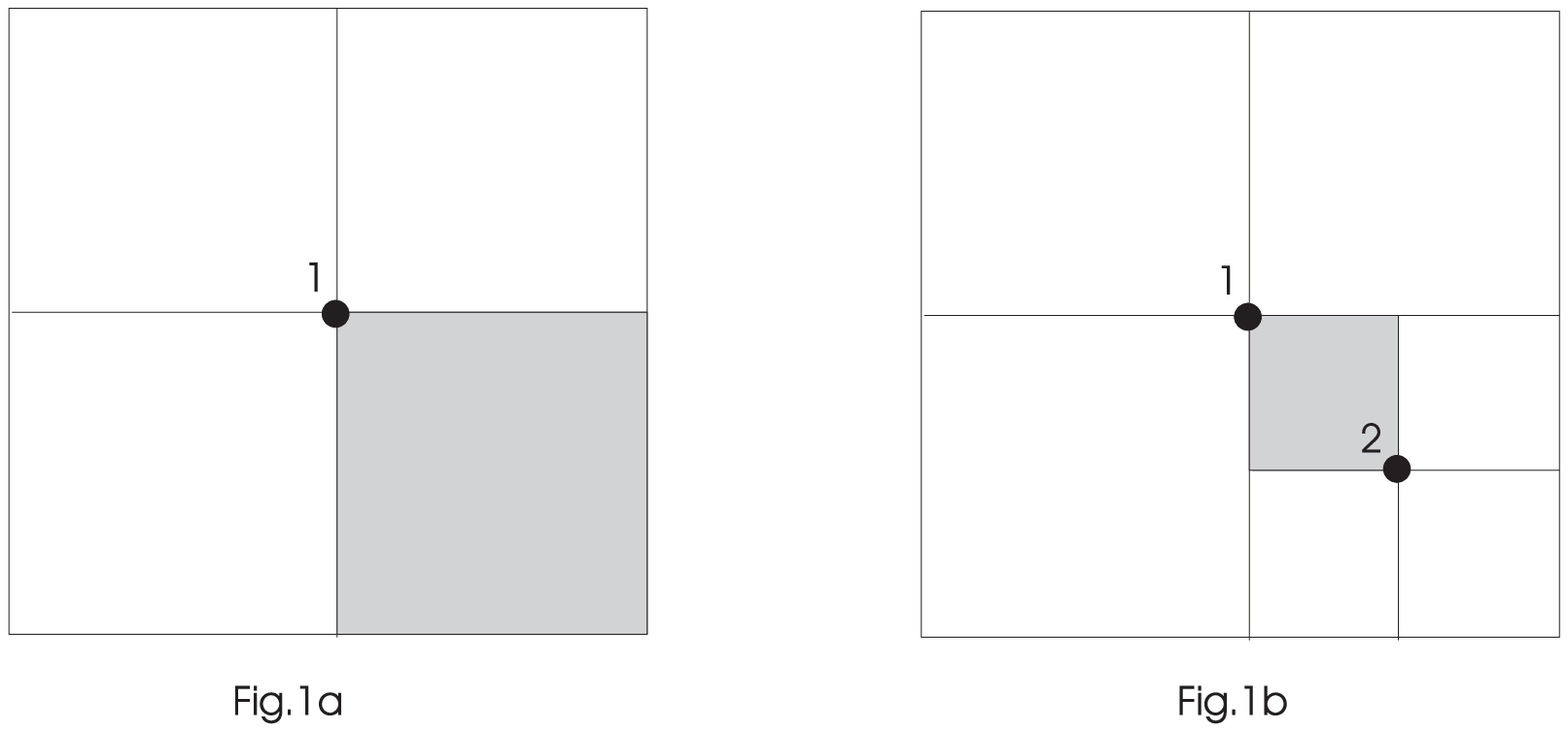,width=0.80\columnwidth,silent=yes}}
 \caption{  }
\end{figure}

\begin{figure}[b]
\vspace*{-5cm}\centerline{\psfig{figure=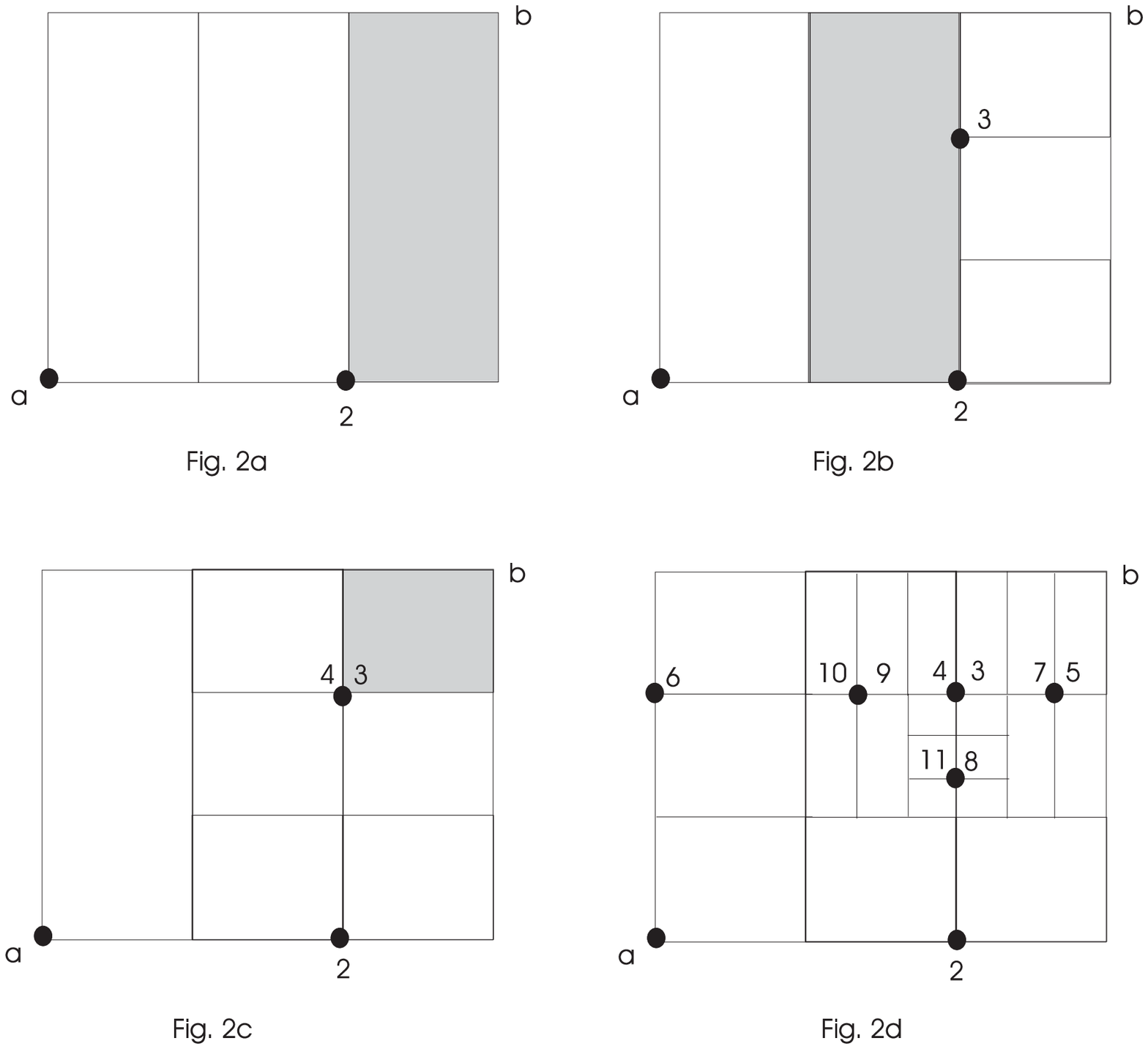,width=0.80\columnwidth,silent=yes}}
 \caption{ }
\end{figure}

\end{document}